\documentclass[11pt]{article}
\usepackage[margin=4.46cm]{geometry}
\usepackage{graphicx}
\usepackage{epsfig}
\usepackage{caption}
\usepackage{amsmath,amsthm,amssymb,dsfont,mathtools,bbm,bm}
\usepackage[T1]{fontenc}
\usepackage{cases}
\usepackage{enumerate}
\usepackage{graphicx,color}
\usepackage{anyfontsize}
\usepackage{float}
\newtheorem{theorem}{Theorem}
\newtheorem{prop}{Proposition}[theorem]
\newtheorem{corollary}{Corollary}[theorem]
\theoremstyle{definition}
\newtheorem{exmp}{Example}[section]
\newcommand{\ds}{\displaystyle}
\newcommand{\mN}{\mathbb{N}}
\newcommand{\mZ}{\mathbb{Z}}

\newcommand{\mP}{\mathbb{P}}
\newcommand{\mE}{\mathbb{E}}

\newcommand{\mF}{\mathcal{F}}

\newcommand{\be}{\begin{equation}}
\newcommand{\ee}{\end{equation}}
\newcommand{\ba}{\begin{eqnarray}}
\newcommand{\ea}{\end{eqnarray}}
\newcommand{\bi}{\begin{itemize}}
\newcommand{\ei}{\end{itemize}}

\newcommand{\real}{\mathbb{R}}

\begin{document}
\title{Enhanced Diffusivity in Perturbed Senile
	Reinforced Random Walk Models}
\author{Thu Dinh 
        and Jack Xin
\thanks{
The authors were partially supported by NSF grants DMS-1522383 and IIS-1632935. 
They are with the Department of Mathematics, University of California, Irvine, CA, 92697, USA.
E-mail: thud2@uci.edu; jack.xin@uci.edu.
}}
\date{}
\maketitle

\begin{abstract}
We consider diffusivity of random walks with transition probabilities depending on 
the number of consecutive traversals of the last traversed edge, the so called senile 
reinforced random walk (SeRW). In one dimension,  
the walk is known to be sub-diffusive with identity reinforcement function. 
We perturb the model by introducing a small probability $\delta$ of 
escaping the last traversed edge at each step. 
The perturbed SeRW model is diffusive for any $\delta >0 $, with 
enhanced diffusivity ($\gg O(\delta^{2})$) in the small $\delta$ regime. We further study stochastically perturbed 
SeRW models by having the last edge escape probability of the form $\delta\, \xi_n$ with $\xi_n$'s being 
independent random variables. 
Enhanced diffusivity in such models are logarithmically close to the so called residual diffusivity 
(positive in the zero $\delta$ limit), with diffusivity  
between $O\left(\frac{1}{|\log\delta |}\right)$ and $O\left(\frac{1}{\log|\log\delta|}\right)$. 
Finally, we generalize our results to higher dimensions where the unperturbed model is already diffusive. 
The enhanced diffusivity can be as much as $O(\log^{-2}\delta)$.
\end{abstract}

\hspace{.12 in} {\bf Key Words:} 
Reinforced random walk, symmetric perturbation,

\hspace{.12 in} enhanced diffusivity, asymptotic analysis.

\vspace{.2 in}

\hspace{.12 in} {\bf AMS Subject Classification:} 60G50, 60H30, 58J37.
\thispagestyle{empty}

\section{Introduction}\label{intro}
\setcounter{equation}{0}
\setcounter{page}{1}
Enhanced diffusivity arises in large scale 
fluid transport through chaotic and turbulent flows, and 
has been studied for nearly a century, see \cite{T21,Kr70,KP79,BCVV95,FP96,MK99,Mur17,LXY17,LXY19} among others. 
It refers to the much larger macroscopic effective diffusivity ($D^E$) than 
the microscopic molecular diffusivity ($D_0$) as the latter approaches zero. 
An example of smooth chaotic flow is  
the time periodic Hamiltonian flow ($X=(x,y) \in \real^2$):
\be\label{tdcell}
\boldsymbol{v}(X,t)=(\cos (y),\cos (x) ) + \theta\; \cos (t)\;(\sin (y),\sin (x)),
\quad
\theta \in (0,1].
\ee 
The first term of (\ref{tdcell}) is a steady flow consisting of periodic arrays of counter-rotating vortices, 
and the second term is a time periodic perturbation that injects an 
increasing amount of disorder into the flow trajectories as $\theta$ becomes larger. 
At $\theta =1$, the flow is fully mixing, and empirically sub-diffusive \cite{ZCX15}. 
The flow (\ref{tdcell}) is one of the simplest models of 
chaotic advection in Rayleigh-B\'enard experiment \cite{CW91}. 
The motion of a diffusing particle in the flow (\ref{tdcell}) 
satisfies the stochastic differential equation (SDE):
\be \label{sde1}
dX_t = \boldsymbol{v}(X_t,t)\, dt + \sqrt{2\,D_0}\, dW_t,\;\; X(0)=(x_0,y_0) \in \real^2, 
\ee
where $W_t$ is the standard 2-dimensional Wiener process. The mean square displacement in 
the unit direction $e$ at large times is given by \cite{BLP2011}:
\be \label{sde2}
\lim_{t \uparrow +\infty} \, E(|(X(t) - X(0))\cdot e |^{2})/t = D^E,  
\ee
where $D^E = D^E(D_0,e,\theta) > D_0$ is the effective diffusivity. 
Numerical simulations \cite{BCVV95,Mur17,LXY17,LXY19} based 
on the associated Fokker-Planck equations (or cell problems of homogenization \cite{BLP2011}) 
suggest that at $e=(1,0)$, $\theta = 1$, $D^E = O(1)$ as $D_0 \downarrow 0$, 
the {\it residual diffusivity} emerges. In fact, $D^E = O(1)$ for $e=(0,1)$ and 
a range of values in $\theta \in (0,1)$ as 
well \cite{LXY17,LXY19}. Recently, computation of (\ref{sde1})-(\ref{sde2}) by structure 
preserving schemes \cite{WXZ_18} reveals residual diffusivity also for a 
time stochastic version of (\ref{tdcell}). At $\theta=0$, enhanced 
$D^{E}$ scales as $O(\sqrt{D_0})\gg D_0$ as $D_0 \downarrow 0$, see \cite{FP94,H03,NPR05} for various 
proofs and generalizations. 
\medskip
 
Motivated by enhanced diffusion in advecting fluids, we are interested in the enhanced diffusion phenomenon in discrete stochastic dynamics such as 
random walk models with some memory or tendency to return. The memory effects on a walker induce a slowdown of transport (movement)  
similar to spinning vortices in fluid flows.  
We shall add a small probability of symmetric random walk and examine the large time behavior of 
the second moment, in similar spirit to (\ref{sde2}). The first work along this line of inquiry
is \cite{LXY_rd} where the baseline model is the so called 
elephant random walk model with stops (ERWS) \cite{ST04,NUK10}. The ERWS  
is non-Markovian and exhibits sub-diffusive, diffusive and super-diffusive regimes. 
The ERWS plays the role of 
flow (\ref{tdcell}). A transition from sub-diffusive to enhanced diffusive regime emerges 
with diffusivity strictly above that of the baseline model (hence residual diffusivity appears) 
as the added probability of symmetric random walk tends to zero \cite{LXY_rd}.   
\medskip

In this paper, we study enhanced diffusivity by perturbing the so called nearest-neighbor reinforced 
senile random walk model (SeRW, \cite{Holmes}) on $\mZ^d$. The model involves a reinforcement 
function $f: \mN \to [-1,\infty)$. The walk $\{S_n\}_{n \geq 0}$ starts at the origin and 
initially steps to one of the $2\, d$ nearest neighbors with equal probability. 
Subsequent steps are defined by the number of times the current undirected edge 
has been traversed consecutively: If $\{S_{n-1},S_n\}$ has been traversed $m$ consecutive times 
in the immediate past, then the probability of traversing that edge in the next step 
is $\frac{1+f(m)}{2d+f(m)}$, with the rest of the possible $2d-1$ choices being equally likely. 
As soon as a new edge is traversed, the reinforcement ends on the previous edge and 
restarts on the new edge. For identity reinforcement function $f$, the walk is sub-diffusive in $d=1$, 
and diffusive in higher dimension \cite{Holmes}. Our work analyzes the asymptotics of the enhanced diffusivity  
when adding a variety of symmetric random walks at small probability.   
\medskip

The rest of the paper is organized as follows. In section 2, we review the baseline SeRW model and 
the key results of \cite{Holmes}. In section 3, we introduce the perturbed SeRW models, 
in which the walk becomes diffusive. In section 4, we state and discuss our main results on  
the diffusivity of the random walk in the perturbed models and the corresponding asymptotics 
for both $d=1$ and $d \geq 2$. The enhancements come logarithmically close to residual diffusivity. 
In section 5, we present proofs of the main results. Concluding remarks are in section 6. 
\medskip

\section{Nearest Neighbor SeRW Model}
A {\it nearest-neighbor senile reinforced random walk} in $\mZ^d$ is a 
sequence $\{S_n\}_{n \geq 0}$ of $\mZ^d$-valued random variables on a 
probability space $(\Omega, \mathcal{F}, \mP_f)$, with corresponding 
filtration $\{\mF_n = \sigma(S_0,...,S_n)\}_{n \geq 0}$, defined by:\\
	$\text{}\quad\bullet$ The walk begins at the origin of $\mZ^d$, i.e. $S_0=0, \mP_f$-almost surely.\\
	$\text{}\quad\bullet$ $\mP_f(S_1=x) = D(x)$, where $D(x) = (2d)^{-1}\mathbbm{1}_{|x|=1}$\\
	$\text{}\quad\bullet$ For $n \in \mN, e_n = \{S_{n-1},S_n\}$ is an $\mF_n$-measurable {\it undirected} edge and
	\[ m_n = \max\{k \geq 1: e_{n-l+1}=e_n \text{ for all } 1 \leq l \leq k \} \]
	is an $\mN$-valued, $\mF_n$-measurable random variable.\\
	$\text{}\quad\bullet$ For $n \in \mN$ and $x \in \mZ^d$ such that $|x|=1$:
	\[
	\mP_f(S_{n+1}=S_n+x \vert \mF_n) =
	\begin{cases}
	\frac{1+f(m_n)}{2d + f(m_n)}, \quad \text{if }\{S_n,S_n+x\} = e_n,\\
	\frac{1}{2d + f(m_n)}, \quad \text{if }\{S_n,S_n+x\} \ne e_n,\\
	\end{cases}\]	

 We shall consider the case $f(m_n) = m_n$, and suppress 
    the $f$ dependence in the probability $\mP_f$ notation. We shall refer to the analysis of SeRW model by 
    Holmes and Sakai \cite{Holmes} and their main results without proofs.\\
	
Let $\tau = \sup\{n \geq 1: S_m = 0 \text{ or } S_1 \text{ } \forall m \leq n\}$ denote 
    the number of times that the walk traverses the first edge before leaving that edge 
  for the first time. Note that $\tau$ is not a stopping time 
   (however $\tau+1 = \inf\{n \geq 2: S_n \ne S_{n-2}\}$ is a stopping time).
	Let $N_x$ denote the number of times the walk $S_n$ visits $x$. 
  If $\mP(N_x=\infty)=1$ for all $x$, we say the walk is recurrent (I). 
 If $\mP(N_x=\infty)=0$ for all $x$, we say the walk is transient (I). 
 If $\mE[N_x] = \infty$ for every $x$, we say the walk is recurrent (II), 
  and if $\mE[N_x] < \infty$ for all $x$, we say the walk is transient (II). 
  Note that for the standard random walk, the two characterizations of recurrence/transience are equivalent; 
 and the walk is recurrent in $d \leq 2$, and transient otherwise. 
   For the senile reinforced random walks, the two notions need not be the same.
	\begin{theorem} (Holmes and Sakai \cite{Holmes})\label{Recurrence}
		For $f$ satisfying $\mP_f(\tau = \infty)=0$, but excluding the degenerate case where $d=1$ and $f(1)=-1$, we have:\\
		(1) $SeRW_{f}$ is recurrent (I)/transient (I) if and only if $SeRW_0$ is recurrent (I)/transient (I).\\
		(2) When $\mE_f[\tau] < \infty, SeRW_f$ is recurrent (II)/transient (II) if and only if $SeRW_0$ is recurrent (II)/transient (II).\\
		(3) When $\mE_f[\tau] = \infty, SeRW_f$ is recurrent (II).
	\end{theorem}
	A consequence of this proposition is the following corollary:
	\begin{corollary}\label{RecurrenceCor}
		The nearest-neighbor senile reinforced random walk with linear reinforcement of the 
      form $f(m) = C\, m$ is recurrent (I), (II) when $d=1,2$ and transient (I) when $d>2$. 
      It is transient (II) for $d > 2$ if and only if $C < 2d-1$.
	\end{corollary}		
	The diffusion constant is defined as $\nu = \lim_{n \to \infty} \mE[|S_n|^2]$ (=1 for the standard 
     random walk) whenever this limit exists. The main result of \cite{Holmes} is:
	\begin{theorem}(Holmes and Sakai \cite{Holmes})\label{DiffusionConstant}
		Suppose that there exists $\epsilon>0$ and $\mE[\tau^{1+\epsilon}]<\infty$. 
        Then the walk is diffusive and the diffusion constant is given by 
		\begin{equation}\label{DiffusionConstantGeneralFormula}
		\nu = \frac{\mP(\tau \text{ odd})}{1-\frac{1}{d}\mP(\tau \text{ odd})}\frac{1}{\mE[\tau]}.
		\end{equation}
	\end{theorem}
	The proof of Theorem \ref{DiffusionConstant} is based on the formula for the Green's function, 
   and a Tauberian theorem, whose application requires the $(1+\epsilon)$th moment of $\tau$ to be finite. 
   Except for the degenerate case, it was shown in \cite{Holmes} that 
   the result holds for all $f$ by a time-change argument. When $\mE[\tau] = \infty$, 
  the right-hand side of \eqref{DiffusionConstantGeneralFormula} is zero, 
  which suggest that the walk is sub-diffusive. \\
	
	When $f(m)=m$, special hypergeometric functions are applicable and 
      various well-known properties of these functions enable a proof of:
	\begin{prop}(Holmes and Sakai \cite{Holmes})\label{NoDiffusion}
		The diffusion constant $\nu$ of the nearest-neighbor senile random walk with 
        reinforcement $f(l)=l$ satisfies $0 < \nu < 1$ when $d > 1$. 
     For the one-dimensional nearest-neighbor model,
		\[ \lim_{n \to \infty} \frac{\log n}{n} \mE[|S_n|^2] = \frac{1-\log 2}{2\log2 - 1}. \]
	\end{prop}
	Hence at $d=1$, the walk is sub-diffusive, slower than diffusion by a logarithmic factor $(\log n)^{-1/2}$. 

\section{Perturbed SeRW Models}
	\subsection{Deterministic Perturbation (Model I)}
	The one-dimensional model with $f(m)=m$ is sub-diffusive. This is partly due to the 
walk having a strong tendency to return to the last traversed edge. 
We add a small perturbation $\delta$ to the conditional probability of $S_{n+1}$ as: 
	\[
	\mP(S_{n+1}=S_n+x \vert \mF_n) =
	\begin{cases}
	\frac{1+m_n}{2 + m_n}-\delta, \quad \text{if }\{S_n,S_n+x\} = e_n,\\
	\frac{1}{2 + m_n}+\delta, \quad \text{if }\{S_n,S_n+x\} \ne e_n.\\
	\end{cases}\]
	In other words, at each step we add a small probability $\delta$ of escaping the last traversed edge, 
where $\delta > 0$ is deterministic. As $m_n \to \infty, \frac{1}{2+m_n} \to 0$. 
So if an edge has already been traversed consecutively too many times, 
the probability of escaping will be dominantly determined by $\delta$. 
This means that the perturbed model will gradually converge to a simplified model where the probability of returning to the last traversed edge is $1-\delta$.
%	It will be shown that the walk is diffusive under this model. 
%However, it is not strong enough to sustain residual diffusivity as $\delta$ tends to zero.
		\subsection{Stochastic Perturbation}	
	\subsubsection{Sequence of i.i.d. perturbations (Model II)}
	%In hope of producing residual diffusivity, we introduce a new model with stochastic perturbation. 
Let $(\xi_n)_{n \in \mN}$ be a sequence of independent identically distributed (i.i.d.) non-negative random variable and consider: 
	\[
	\mP(S_{n+1}=S_n+x \vert \mF_n) =
	\begin{cases}
	\max\{\frac{1+n}{2 + n}-\delta\xi_n,0\}, \quad \text{if }\{S_n,S_n+x\} = e_n,\\
	\min\{\frac{1}{2 + n}+\delta\xi_n,1\}, \quad \text{if }\{S_n,S_n+x\} \ne e_n.\\
	\end{cases}
	\]
	At each step, the random variable $\xi_n$ takes a value, then the reinforcement is based on this value. We only assume that $\xi_n$ is continuous with probability density function $f = f_{\xi_n}$.\\
	
 Notice that if $\xi_n$ takes any value greater than $\frac{1+n}{2+n}$,  the walk will escape the 
last traversed edge on the $n+1^{th}$ turn. So in this model, the tail of the distribution function $f$ 
provides a stronger chance of breaking out of the last traversed step, leading to more enhanced diffusion.
	\subsubsection{Sequence of independent perturbations (Model III)}
	To further enhance diffusivity, we shall consider the situation that $(\xi_n)_{n \in \mN}$ are 
no longer i.i.d., but rather have $n$-dependent distributions.
	\[
	\mP(S_{n+1}=S_n+x \vert \mF_n) =
	\begin{cases}
	\max\{\frac{1+n}{2 + n}-\delta\xi_n,0\}, \quad \text{if }\{S_n,S_n+x\} = e_n,\\
	\min\{\frac{1}{2 + n}+\delta\xi_n,1\}, \quad \text{if }\{S_n,S_n+x\} \ne e_n,\\
	\end{cases}
	\]
	For example, $\xi_n$'s can have the same type of distribution and expectation, but with variance $n^2$. This modification will reinforce the probability of the walk breaking out of the last traversed edge. We only assume that $\mE[\xi_n] < \infty$, for all $n$.

	\section{Main Results}
	The diffusivity from the perturbation (the simple symmetric random walk) similar to "molecular diffusivity'' 
$D_0$ of (\ref{sde1}) is $\nu_\delta = \delta^2$. We will show that, 
in all of our three models, the enhanced diffusivity is much greater than $O(\delta^2)$. 
Our main results are stated in the following theorems.
	\begin{theorem}\label{DeterministicResult}
		The deterministic perturbed model (I) is diffusive for any $\delta > 0$, and the diffusion constant is given by 
		\begin{equation}\label{NuFormula}
		\nu = \frac{\mP(\tau \text{ odd})}{\mP(\tau \text{ even})\mE[\tau]}.
		\end{equation}
		Moreover,
		\begin{equation}
		\nu(\delta) = O\left(\frac{1}{|\log\delta |}\right) \quad \text{as} \quad \delta \to 0^+.
		\end{equation}
		%Hence the diffusion is not strong enough to produce residual diffusivity (although $\nu$ is bounded away from zero) as $\delta$ tends to zero.
	\end{theorem}
	The formula (\ref{NuFormula}) for $\nu$ is a direct result of Theorem \ref{DiffusionConstant}. 
It is dramatic that the walk becomes diffusive for any value of $\delta > 0$. 
Proposition \ref{NoDiffusion} says the walk is sub-diffusive by an order of $\log n$. 
The added perturbation reduces the probability of revisiting the last traversed edge. 
However small, the perturbation is enough to create diffusivity. 
%However, this model does not produce any residual diffusivity as $\delta$ tends to zero.\\

	To prove Theorem \ref{DeterministicResult}, we first verify that the model is diffusive by checking the condition of Theorem \ref{DiffusionConstant}, then we will find a lower bound for $\mE[\tau]$ and show that the bound goes to $\infty$ as $n \to \infty$. A straightforward computation shows $1 \leq \frac{\mP(\tau \text{ odd})}{\mP(\tau \text{ even})} \leq 2$. This concludes the proof. In the last section, we discuss the rate at which $\nu$ goes to zero as $\delta$ tends to zero.
	
	\begin{theorem}\label{StochasticIIDResult}
		The stochastic perturbed model (II) is diffusive for any $\delta > 0$, and 
the diffusion constant is given by the same formula as in Theorem \ref{DeterministicResult}. Moreover,\\ 
%However, the model is not strong enough to produce residual diffusivity as $\delta$ tends to zero.\\
		(i) If $\mE[\xi_n] < \infty$, then  $\nu(\delta) = O\left(\frac{1}{|\log\delta|}\right)$ as $\delta \to 0^+$.\\
		(ii) If $\mE[\xi_n] = \infty$, one can construct  $\xi_n$ so that $\nu(\delta) = O(\frac{1}{\log|\log\delta|})$ as $\delta \to 0^+$.
	\end{theorem}
	
	Similar to the deterministic case, the stochastic perturbed model is still not strong enough 
to sustain residual diffusivity. We can, however, reduce the rate at which $\nu$ converges to 0. 
If $\xi_n$ has infinite expected value (fat tail), then $\xi_n$ is more likely to attain very large values, 
and the walk is less likely to get stuck. The maximal enhancement on $\nu(\delta)$ is $O(\frac{1}{\log|\log\delta|})$.

	\begin{theorem}\label{StochasticVarResult}
		The stochastic perturbed model (III) is diffusive for any $\delta> 0$.  
The diffusion constant is given by the same formula as in Theorem \ref{DeterministicResult} 
 %However, the model is not strong enough to produce residual diffusivity, and 
with $\nu(\delta) = O\left(\frac{1}{|\log\delta |}\right)$ as $\delta \to 0^+$.
	\end{theorem}
	
	The proofs of the three theorems above are based on Theorem \ref{DiffusionConstant} to show diffusivity and  
the calculation of the diffusion constant $\nu$. Our approach is elementary and relies heavily on 
the computation of the quantity $\mP(\tau \geq n)$. The absence of residual diffusivity 
and the rate of convergence are obtained via asymptotic analysis in the small $\delta$ regime.
	
	\begin{theorem}\label{HigherDimensionResult}
		When the baseline diffusive SeRW model on $\mZ^d$ ($d \geq 2$) is perturbed into models (I, II, III), 
        we have the following:\\
		(i) Under model I, the walk has a linearly enhanced diffusivity: 
		\[ \nu_\delta = \nu_0 + O(\delta), \]
       where $\nu_0$ is the diffusivity of the unpeturbed model.\\
	%	where $M := \mE[\tau]$ of the unperturbed model.\\
		(ii) Under models II and III, if $\mE[\xi_n] < \infty$, for all $n$, the walk has the same linear enhanced diffusivity as in model I. \\
		(iii) Under models II and III, if $\mE[\xi_n] = \infty$, for all $n$, 
         one can construct $\xi_n$ to achieve the following enhanced diffusivity rates:
		\begin{align*}
		(a) &\quad \nu_\delta = \nu_0 +  O(\delta\, |\log\delta |),\\
		(b) &\quad \nu_\delta = \nu_0 + O(\delta^{j}), \quad\text{for some }\;\; j \in (0,1),\\
		(c) &\quad \nu_\delta = \nu_0 + O(\log^{-2}\delta).
		\end{align*}
		
	\end{theorem}
	
	\section{Proofs of Main Results}
	\subsection{Theorem \ref{DeterministicResult}: Existence of positive diffusion constant}\label{DeterministicProof}
	First we verify the perturbed model is diffusive. It is easy to see that
	\[ \mP(\tau=1) = \frac{1}{3}+\delta \quad\quad\text{and}\quad\quad \mP(\tau=n) = 
\left[\prod_{k=2}^n\left(\frac{k}{k+1}-\delta\right)\right]\left( \frac{1}{n+2}+\delta \right) \]
	for $n \geq 2$. We will show there exists $\epsilon > 0$ such that $\mE[\tau^{1+\epsilon}] < \infty$ and 
apply Theorem \ref{DiffusionConstant}. The following is an upper bound for $\mP(\tau = n)$ when $n \geq 2$:\\
	\begin{align*}
	\mP(\tau = n) &= \left[ \prod_{k=2}^n\left(\frac{k}{k+1}-\delta\right) \right] \left(\frac{1}{n+2}+\delta\right)\\
	&= \left(\frac{2}{3}-\delta\right)\left(\frac{3}{4}-\delta\right)...\left(\frac{n}{n+1}-\delta\right)\left(\frac{1}{n+2}+\delta\right)\\
	&= \frac{2(1-\frac{3\delta}{2})3(1-\frac{4\delta}{3})...n(1-\frac{(n+1)\delta}{n})}{3\cdot 4 ... \cdot (n+1)}\left(\frac{1}{n+2}+\delta\right)\\
	&\leq \frac{2}{n+1} e^{-\frac{3\delta}{2}}e^{-\frac{4\delta}{3}}...e^{-\frac{(n+1)\delta}{n}}\left(\frac{1}{n+2}+\delta\right)\\
	&= \frac{2}{n+1} \exp\left\{ -\sum_{k=2}^n \delta\left(1+\frac{1}{k}\right)\right\}\left(\frac{1}{n+2}+\delta\right)\\
	&\leq \frac{2}{n+1} \exp\left\{ \delta( -n + 1 - \log n + 1) \right\}\left(\frac{1}{n+2}+\delta\right)\\
	&= \frac{2e^{2\delta}(1+(n+2)\delta)}{(n+1)(n+2)e^{\delta n}n^\delta}
	\end{align*}
	where the first inequality follows since $1-x \leq e^{-x}$ for all $x$, and the second inequality since $\log n \leq \sum_{k=1}^n \frac{1}{n}$. Letting $\epsilon = \delta$, we have
	\begin{align*}
	\mE[\tau^{1+\delta}] &= \sum_{n=1}^\infty n^{1+\delta}\mP(\tau=n)\\
	&= \frac{1}{3}+\delta + \sum_{n=2}^\infty \frac{2e^{2\delta}n(1+(n+2)\delta)}{e^{\delta n}(n+1)(n+2)} < \infty
	\end{align*}
	Thus by Theorem \ref{DiffusionConstant}, the walk is diffusive.\\
%	We wish to obtain a residual diffusion constant as $\delta$ tends to zero. 
%However, the model is not strong enough to sustain the diffusivity. 
In the second part of this proof, we will show $\nu \to 0$ as $\delta \to 0^+$. By Theorem \ref{DiffusionConstant}, the diffusion constant simplifies to
	\[ \nu = \frac{\mP(\tau \text{ odd})}{\mP(\tau \text{ even})\mE[\tau]}. \]
	It suffices to show $\mE[\tau] \to \infty$ as $\delta \to 0^+$. To that end, it is more convenient to use the formula $\mE[\tau] = \sum_{n=1}^\infty \mP(\tau \geq n)$. We have
	\[ \mP(\tau \geq 1) = 1 \quad\quad\text{and}\quad\quad \mP(\tau \geq n) = \prod_{k=2}^n\left(\frac{k}{k+1}-\delta\right) \]
	for $n \geq 2$. The following computation gives a lower bound for $\mP(\tau \geq n)$ when $n \geq 2$:
	\begin{align*}
	\mP(\tau \geq n) &= \prod_{k=2}^n\left(\frac{k}{k+1}-\delta\right)\\
	&= \frac{2(1-\frac{3\delta}{2})3(1-\frac{4\delta}{3})...n(1-\frac{(n+1)\delta}{n})}{3\cdot 4 ... \cdot (n+1)}\\
	&\geq \frac{2}{n+1} e^{-2(\frac{3\delta}{2})}e^{-2(\frac{4\delta}{3})}...e^{-2(\frac{(n+1)\delta}{n})}\\
	&= \frac{2}{n+1} \exp\left\{ -2\sum_{k=2}^n \delta\left(1+\frac{1}{k}\right)\right\}\\
	&\geq \frac{2}{n+1} \exp\left\{ -2\delta( n -2 + \log n + \gamma) \right\}\\
	&= \frac{2e^{4\delta}}{(n+1)e^{2\delta\gamma}e^{2\delta n} n^{2\delta}}\\
	&\geq \frac{2e^{4\delta}}{2ne^{2\delta\gamma}e^{2\delta n} n^{2\delta}}
	\end{align*}
	where the first inequality follows since $1-x \geq e^{-2x}$ holds for small $x \geq 0$, and the second equality since $\sum_{k=1}^n \frac{1}{k} \leq \log n + \gamma$, where $\gamma$ is the Euler constant. \\
	It remains to show $\sum_{n=1}^\infty \frac{2e^{4\delta}}{(n+1)e^{2\delta\gamma}e^{2\delta n} n^{2\delta}} \to \infty$ as $\delta \to 0^+$. Since the terms in the summation are positive and decreasing, we can use the integral test for convergence. After multiplying by a constant, 
it suffices to compute
	\[ \int_1^\infty \frac{e^{-2\delta x}}{x^{1+2\delta}}dx. \]
	Letting $t = -2\delta$, we have
	\[ \int_1^\infty \frac{e^{-2\delta x}}{x^{1+2\delta}}dx
	= \int_{2\delta}^\infty \frac{e^{-t}}{\left(\frac{t}{2\delta}\right)^{1+\delta}} \frac{dt}{2\delta}
	= (2\delta)^{\delta} \int_{2\delta}^\infty \frac{e^{-t}}{t^{1+2\delta}}dt
	= (2\delta)^{\delta} \Gamma(-2\delta,2\delta) \]
	where $\Gamma(\cdot,\cdot)$ is the Incomplete Upper Gamma function \cite{Milton}. It is straightforward to 
      verify that $(2\delta)^\delta \to 1$ as $\delta \to 0^+$. By \cite{Milton}, 
     $\Gamma(-2\delta,2\delta) \to \infty$ as $\delta \to 0^+$.\\
	Thus, we have shown that a lower bound for $\mE[\tau]$ diverges as $\delta$ tends to 0. 
    By Theorem \ref{DiffusionConstant}, $\nu$ converges to 0. 
    Therefore the perturbed model is not strong enough to sustain a residual diffusivity.
	\subsection{Rate of convergence}\label{RateOfDivergence}
	Since a residual diffusion is not achievable, it is natural to ask how fast $\nu$ is decreasing as $\delta$ tends to 0. In this section, we will verify that in the perturbed model, the diffusivity converges to 0 
at a rate of $\frac{1}{|\log \delta|}$.\\
	Let $k = 2\delta$ and consider the integral above as a function of $k$, i.e.,
	\begin{equation}\label{RateFucntion}
	f(k) = \int_1^\infty \frac{e^{-k x}}{x^{1+k}} dx
	\end{equation}
	Then
	\[ f'(k) = -\int_1^\infty \frac{e^{-kx}}{x^{1+k}}(x+\log x) dx \]
	Since $x \gg \log x$ as $x \to \infty$, $f'(k)$ is dominantly determined by the term with $x$, namely
	\[
	f'(k) \sim -\int_1^\infty \frac{e^{-kx}}{x^k}dx
	\]
	let $u = x^{1-k}$, so $du = (1-k)x^{-k}dx$, we have
	\[ f'(k) \sim -\frac{1}{1-k} \int_1^\infty e^{-ku^{\frac{1}{1-k}}}du 
	= -\frac{1}{1-k} \int_1^\infty e^{-(k^{1-k}u)^{\frac{1}{1-k}}}du \]
	let $v = k^{1-k}u$, the integral becomes
	\[ f'(k) \sim -\frac{k^{-1+k}}{1-k} \int_{k^{1-k}}^\infty e^{-v^{\frac{1}{1-k}}}dv \]
	as $k \to 0^+$,
	\[ \int_{k^{1-k}}^\infty e^{-v^{\frac{1}{1-k}}}dv \to \int_0^\infty e^{-v}dv = 1 \]
	thus $f'(k) \sim -\frac{k^{-1+k}}{1-k}$ as $k \to 0^+$. Finally,
	\[ \lim_{k \to 0^+} \frac{-f(k)}{\log(\delta)} = \lim_{k \to 0^+} \frac{-f(k)}{\log k-\log 2}
	\stackrel{\mathclap{\text{L'H}}}{=} \lim_{k \to 0^+} \frac{-f'(k)}{1/k}
	= \lim_{k \to 0^+} \frac{k^{-1+k}k}{1-k} = 1 \]
	An identical computation shows $\lim_{k \to 0^+} \frac{f(\delta)}{-\log\delta} = 1$.
	Since 
	\[ \sum_{n=1}^\infty \frac{2e^{4\delta}}{(n+1)e^{2\delta\gamma}e^{2\delta n} n^{2\delta}} \leq \mE[\tau]
	\leq \sum_{n=1}^\infty \frac{2e^{2\delta}}{(n+1)e^{\delta\gamma}e^{\delta n} n^{\delta}}, \] 
	after multiplying by a constant, we have $\mE[\tau] \sim C_1\, |\log \delta |$. 
Applying the formula of Theorem \ref{DiffusionConstant}, we have $\nu_\delta = O\left(\frac{1}{|\log \delta|}\right)$.
	
	\subsection{Theorem \ref{StochasticIIDResult}: Existence of positive diffusion constant}\label{StochasticProof}
	The formula for the diffusion constant $\nu$ follows directly from Theorem \ref{DiffusionConstant}. The proof of Theorem \ref{DiffusionConstant} is based on the formula for the Green's function, and a standard Tauberian theorem. It utilized the following functions and quantities:
	\[ G_z(x) = \sum_{n=0}^\infty z^n \mP(S_n=x), \quad \text{for } z \in [0,1] \]
	\[\begin{cases}
	a_z = \ds\sum_{n=2}^\infty z^n \mP(\tau \geq n) \mathbbm{1}_{\{n \text{ even}\}}\\
	b_z = \ds\sum_{n=2}^\infty z^n \mP(\tau \geq n) \mathbbm{1}_{\{n \text{ odd}\}}
	\end{cases}
	\quad
	\begin{cases}
	p_z = \ds\sum_{n=1}^\infty z^n \mP(\tau = n) \mathbbm{1}_{\{n \text{ even}\}}\\
	q_z = \ds\sum_{n=1}^\infty z^n \mP(\tau = n) \mathbbm{1}_{\{n \text{ odd}\}}
	\end{cases}\]
	and other variables built up from $a_z,b_z,p_z$, and $q_z$. We will show below that, even though the model is stochastic, $\mP(\tau \geq n)$ is still deterministic. Thus the proof of Theorem \ref{DiffusionConstant} still applies and gives the formula for $\nu$.

	Given that an edge has been traversed $n$ times, let $P_n$ denote the total probability of breaking out of this edge on the $(n+1)^{th}$ turn, and let $Q_n$ denote the probability of traversing this edge again on the $(n+1)^{th}$ turn. Then $P_n$ is the sum of all the terms of the form $\frac{n+1}{n+2}-\delta\xi$, given that $\xi_n=\xi \leq \frac{n+1}{\delta(n+2)}$. Formally,
	
	\[ P_n = \int_0^{\frac{n+1}{\delta(n+2)}} \left(\frac{n+1}{n+2}-\delta x\right)f(x)	 dx \]
	and 
	\[ Q_n = \left(\int_0^{\frac{n+1}{\delta(n+2)}} \left(\frac{1}{n+2} - \delta x\right)f(x)dx \right) + \mP\left(\xi_n > \frac{n+1}{\delta(n+2)}\right) \]
	Similar to the previous result, for $n \geq 2$, we have
	\[ \mP(\tau=n) = \left(\prod_{i=1}^{n-1}P_i\right)Q_n \hspace{.6in}
	\text{and} \hspace{.6in} \mP(\tau\geq n) = \prod_{i=1}^{n-1}P_i. \]
	An upper bound for $\mP(\tau \geq n)$ is
	\begin{align*}
	\mP(\tau\geq n) &= \prod_{i=1}^{n-1}\left( \int_0^{\frac{i+1}{\delta(i+2)}} \left(\frac{i+1}{i+2}-\delta x\right)f(x)dx \right)\\
	&\leq \prod_{i=1}^{n-1}\left( \frac{i+1}{i+2} - \delta\int_0^{\frac{i+1}{\delta(i+2)}} xf(x)dx \right)\\
	&\leq \prod_{i=1}^{n-1}\left( \frac{i+1}{i+2} - \delta \int_0^{\frac{2}{3\delta}} xf(x)dx \right).
	\end{align*}
	Let $\mu :=  \delta\int_0^{\frac{2}{3\delta}} xf(x)dx$. Then $\mu$ is a constant for each fixed $\delta$. Thus $\mP(\tau\geq n) = \prod_{i=1}^{n-1}\left( \frac{i+1}{i+2} - \mu \right)$, which has the same form as in the deterministic case. By a similar computation, there exists $\epsilon>0$ such that $\mE[\tau^{1+\epsilon}] < \infty$, and the walk is diffusive.\\

	Recall Theorem \ref{DiffusionConstant}, the diffusion constant is
	\[ \nu = \frac{\mP(\tau \text{ odd})}{\mP(\tau\text{ even})}\frac{1}{\mE[\tau]} \]
	In order to sustain residual diffusivity, we need $\mE[\tau] \not\to \infty$ as $\delta \to 0^+$. Using the formula $\mE[\tau] = \sum_{n=1}^\infty \mP(\tau \geq n)$, we get
	\begin{equation}\label{ETauStochastic}
	\mE[\tau]
	= 1 +  \sum_{n=2}^\infty\prod_{i=1}^{n-1}\left( \int_0^{\frac{i+1}{\delta(i+2)}}
 \left(\frac{i+1}{i+2}-\delta x\right)f(x)dx \right).
	\end{equation}
	Suppose $\mE[\xi_n]<\infty$. Then by Fatou's lemma,
	\begin{align*}
	&\liminf_{\delta \to 0^+} \mE[\tau] = \liminf_{\delta \to 0^+}\left(
	1 +  \sum_{n=2}^\infty\left[\prod_{i=1}^{n-1} \int_0^{\frac{i+1}{\delta(i+2)}} \left(\frac{i+1}{i+2}-\delta x\right)f(x)dx \right]\right)\\
	&\geq 1 +  \sum_{n=2}^\infty \liminf_{\delta \to 0^+}\prod_{i=1}^{n-1} \left[ \left(\frac{i+1}{i+2}\right)\int_0^{\frac{i+1}{\delta(i+2)}}f(x)dx -\delta \int_0^{\frac{i+1}{\delta(i+2)}}xf(x)dx \right]\\
	&= 1 + \sum_{n=2}^\infty \liminf_{\delta \to 0^+}\left[\prod_{i=1}^{n-1} \left(\frac{i+1}{i+2} - \delta \mE[\xi_n]\right)\right]\\
	&= 1 + \sum_{n=2}^\infty \left(\prod_{i=1}^{n-1}\frac{i+1}{i+2}\right)\\
	&= 1 + \sum_{n=2}^\infty \frac{2}{n+1} = \infty.
	\end{align*}
	Since a lower bound for $\mE[\tau]$ diverges to $\infty$, the corresponding 
upper bound for $\nu$ converges to 0. Thus $\nu \to  0$ 
%and the model fails to produce residual diffusivity 
as $\delta \to 0^+$.  Moreover, since $\mE[\xi_n]$ is a finite constant, 
the computation from Section \ref{RateOfDivergence} shows 
$\nu(\delta) = O(\frac{1}{|\log\delta |})$ as $\delta \to 0^+$.

	\subsection{Random variables with infinite expectation}
	\subsubsection{Necessary asymptotic behavior of the pdf of $\bm{\xi_n}$}
	Suppose $\xi_n$ is a random variable with support in $[0,\infty)$ and $\mE[\xi_n] = +\infty$. Let $f = f_{\xi_n}$ be the 
probability density function (pdf) of $\xi_n$, we have
	\[ \int_0^\infty f(x)dx = 1 \quad \text{and} \quad \int_0^\infty xf(x)dx = \infty.  \]
	We will study the asymptotic behavior of such $f$. Since $\int_0^\infty f(x)dx = 1$, we require $f(x) \leq O(x^{-n})$, for some $n>1$.\\
	On the other hand, $\int_0^\infty xf(x)dx = \infty$ implies $xf(x) \geq O(x^{-1})$. Thus, the necessary asymptotic behavior for $f$ is
	\[ O\left(\frac{1}{x^2}\right) \leq f(x) < O\left(\frac{1}{x}\right). \]
	\begin{exmp}\label{CauchyExample}
	A random variable $\xi_n$ with $f(x) = O\left(\frac{1}{x^2}\right)$.\\	
	Let $\xi_n$ be non-negative Cauchy random variables with $x_0=0$ and pdf
	\[ f_{\xi_n}(x) = \frac{2}{\pi\gamma\left[1+\left(\frac{x}{\gamma}\right)^2\right]}
	= \frac{2\gamma}{\pi(x^2+\gamma^2)}. \]
	Then
	\begin{align*}
	&\mP(\tau \geq n)
	=  \prod_{i=1}^{n-1} \left[ \left(\frac{i+1}{i+2}\right)\int_0^{\frac{i+1}{\delta(i+2)}}\frac{2\gamma}{\pi(x^2+\gamma^2)}dx -\delta \int_0^{\frac{i+1}{\delta(i+2)}}\frac{2\gamma x}{\pi(x^2+\gamma^2)}dx \right]\\
	&= \prod_{i=1}^{n-1} \left[ \left(\frac{i+1}{i+2}\right)\int_0^{\frac{i+1}{\delta(i+2)}}\frac{2\gamma}{\pi(x^2+\gamma^2)}dx
	- \delta \left(\frac{\gamma\log(x^2+\gamma^2)}{\pi}\right)\bigg\vert_{x=0}^{x=\frac{i+1}{\delta(i+2)}} \right]\\
	&= \prod_{i=1}^{n-1} \left[ \left(\frac{i+1}{i+2}\right)\int_0^{\frac{i+1}{\delta(i+2)}}\frac{2\gamma}{\pi(x^2+\gamma^2)}dx
	- O\left(\delta\log\frac{1}{\delta}\right)
	\right]
	\end{align*}
	and by Fatou's lemma,
	\begin{align*}
	&\liminf_{\delta \to 0^+} \mE[\tau]
	\geq 1 + \sum_{n=2}^\infty \liminf_{\delta \to 0^+} \mP(\tau \geq n)\\
	&= 1 + \sum_{n=2}^\infty \liminf_{\delta \to 0^+} \prod_{i=1}^{n-1} \left[ \left(\frac{i+1}{i+2}\right)\int_0^{\frac{i+1}{\delta(i+2)}}\frac{2\gamma}{\pi(x^2+\gamma^2)}dx
	- O\left(\delta\log\frac{1}{\delta}\right).
	\right]\\
	&= 1 + \sum_{n=2}^\infty \frac{2}{n+1} = \infty.
	\end{align*}
	Similar to the above result, since a lower bound for $\mE[\tau]$ diverges to $\infty$, we have $\nu \to 0$ as $\delta \to 0^+$. Thus, even though the non-negative Cauchy distribution has a "fat'' tail, the growth rate of $\int_0^{\frac{i+1}{\delta(i+2)}}xf(x)dx$ is still not fast enough to produce residual diffusivity.
	\end{exmp}
	
	\subsubsection{Non-existence of residual diffusivity, rate of convergence}\label{StochasticRate}
	The case  where $f(x) = O(x^{-2})$ was covered in example \ref{CauchyExample}. In general, if
	\[ O\left(\frac{1}{x^2}\right) < f(x) < O\left(\frac{1}{x}\right) \]
	then
	\[ O\left(\frac{1}{x}\right) < xf(x) < O(1) \]
	which implies
	\[ \delta O\left(\log\frac{1}{\delta}\right) < \delta \int_0^{\frac{i+1}{\delta(i+2)}} xf(x)
	< \delta O\left( \frac{1}{\delta} \right). \]
	Taking the limit as $\delta \to 0^+$, we have $\delta \int_0^{\frac{i+1}{\delta(i+2)}} xf(x) \to 0$, which implies
	\[ \mP(\tau \geq n) = \prod_{i=1}^{n-1}\left( \int_0^{\frac{i+1}{\delta(i+2)}} 
\left(\frac{i+1}{i+2}-\delta x\right)f(x)dx \right) \to \frac{1}{n} \text{ as } \delta \to 0^+. \]
	Therefore $\mE[\tau] \to \infty$ and, subsequently, $\nu \to 0$. \\
%Thus the model fails to sustain residual diffusivity in all cases.\\
	
For the asymptotic behavior of $\nu(\delta)$, we study 3 cases:
	\paragraph{Case 1,} $\bm{f(x) = O(x^{-2}):}$\\
	By example \ref{CauchyExample}, as $\delta \to 0^+$, 
	\[ \mP(\tau \geq n) = \prod_{i=1}^{n-1} \left(\frac{i+1}{i+2} - C\delta\log\left(\frac{1}{\delta}\right) \right) \]
	A similar computation to the last part of section \ref{DeterministicProof} shows that, after multiplying by a constant, to compute $\mE[\tau]$, it suffices to compute
	\[ g(\delta\log(1/\delta)) = \int_1^\infty \frac{e^{-\delta\log(1/\delta)}}{x^{1+\delta\log(1/\delta)}}. \]
	And by the computation of section \ref{RateOfDivergence}, which shows $\lim_{\delta \to 0^+} \frac{g(k)}{\log k} = 1$, we have
	\[ \lim_{\delta \to 0} \frac{g(\delta\log(1/\delta))}{\log(\delta\log(1/\delta))} = 1. \]
	This implies 
	\[ \mE[\tau] \sim C_1\log(|\delta\log\delta|) \]
	and therefore
	\[ \nu \sim \frac{C_2}{\log(|\delta\log\delta|)} \sim \frac{C_3}{\log\delta}. \]
	\paragraph{Case 2,} $\bm{f(x) = O(x^{-(1+j)})}$, \textbf{for} $\bm{0 < j < 1:}$\\
	A similar calculation to example \ref{CauchyExample} shows, as $\delta \to 0^+$,
	\[ \mP(\tau \geq n) = \prod_{i=1}^{n-1} \left(\frac{i+1}{i+2} - C\delta^j \right) \]
	and a calculation similar to Case 1 shows
	\[ \mE[\tau] \sim C_1\, |\log(\delta^j)| = C_2\, |\log(\delta)|. \]
	So in this case, 
	\[ \nu \sim \frac{C_3}{|\log(\delta)|} \]
	which is the same result as the deterministic case.
	\paragraph{Case 3,} $\bm{f(x) < O(x^{-(1+j)})}$, \textbf{for any} $\bm{0<j<1}$ \textbf{and} $\bm{f(x) > O(x^{-2}):}$\\
	One such example is $f(x) = O\left(\frac{1}{x(\log x)^2}\right)$. Then
	\[ \int_0^{\frac{i+1}{\delta(i+2)}} xf(x)dx = \int_0^{\frac{i+1}{\delta(i+2)}} \frac{C}{\log^2x}dx \]
	which is a well known logarithm integral with asymptotic behavior:
	\[ \int \frac{1}{\log^2 x}dx = li(x) - \frac{x}{\log x} = O\left(\frac{x}{\log^2 x}\right) \]
	therefore
	\[ \delta\int_0^{\frac{i+1}{\delta(i+2)}} xf(x)dx = \delta O\left( \frac{1}{\delta\log^2\left(\frac{C_1}{\delta}\right)} \right)
	= O\left(\frac{1}{\log^2(\delta)}\right) \]
	as $\delta \to 0^+$. This implies
	\[ \mP(\tau \geq n) = \prod_{i=1}^{n-1} \left(\frac{i+1}{i+2} - \frac{C_2}{\log^2\delta} \right) \]
	and a similar calculation to Case 1 shows
	\[ \mE[\tau] \sim C_3\log(\log^2\delta) = C_4 \log|\log\delta|. \]
	Thus, we have constructed a random variable $\xi_n$ such that $\nu$ converges to zero at a rate of
	\[ \nu \sim \frac{C}{\log|\log\delta|}. \]
	
	\subsection{Proof of Theorem \ref{StochasticVarResult}}
	Theorem \ref{StochasticVarResult} is a consequence of Theorem \ref{StochasticIIDResult}. The fact that the model is diffusive for any $\delta > 0$ follows directly. For the rate at which $\nu$ tends to 0, let $f_n$ be the p.d.f. of $\xi_n$ and recall that
	\[ \mP(\tau\geq n) = \prod_{i=1}^{n-1}\left( \int_0^{\frac{i+1}{\delta(i+2)}} 
\left(\frac{i+1}{i+2}-\delta x\right)f_n(x)dx \right). \]
	Since $\mE[\xi_n] < \infty$ for all $n$, one can find a random variable $Y$ with $\mE[Y]=\infty$ with p.d.f. $f_Y$ such that, for sufficiently small $\delta$,
	\[ \delta\int_0^\frac{i+1}{\delta(i+2)} xf_n(x)dx \leq \delta\int_0^\frac{i+1}{\delta(i+2)} yf_Y(y)dy \]
	so as $\delta \to 0^+$, 
	\[ \mP(\tau \geq n) \geq \prod_{i=1}^{n-1} \left(\frac{i+1}{i+2} - 
\delta\int_0^\frac{i+1}{\delta(i+2)} yf_Y(y)dy \right). \]
	Notice the expression on the RHS matches the case of infinite expectation of the Theorem \ref{StochasticIIDResult}. Therefore $\mE[\tau]$ grows at least as fast as the previous case, and hence so is the decay rate of $\nu_\delta$. One can choose $Y$ so that $f_Y(y) = O(y^{-2})$ (Similar to Case 1 of section \ref{StochasticRate}), so that $\nu_Y(\delta) \sim O(|\log\delta|)$. 
Then $\nu_\delta$ decays at a rate of at most $O(|\log\delta|)$ (by section \ref{StochasticProof}), and at least $O(\log\delta)$, from the previous case. It follows that $\nu_\delta = O(\log\delta)$. 
	
	\subsection{Theorem 6: Results in higher dimensions}
	\subsubsection{Perturbation under model I:}
	For $d \geq 2$, the model becomes 
	\[
	\mP(S_{n+1}=S_n+x \vert \mF_n) =
	\begin{cases}
	\max\{\frac{1+n}{2d + n}-\delta,0\}, \quad \text{if }\{S_n,S_n+x\} = e_n,\\
	\min\{\frac{1}{2d + n}+\delta,1\}, \quad \text{if }\{S_n,S_n+x\} \ne e_n.
	\end{cases}
	\]
	A similar computation to that of the one-dimensional case shows, for $n \geq 2d$,
	\begin{align*}
	& \mP(\tau \geq n) = \prod_{k=2}^n \left(\frac{1+k}{2d+k}-\delta\right)\\
	&= \frac{(2d)!}{(n+2)(n+3)...(n+2d)}
	\left(1-\frac{2d+1}{2}\delta\right)...
	\left(1-\frac{2d+n}{n+1}\delta\right)\\
	&\to \frac{(2d)!}{(n+2)(n+3)...(n+2d)} \exp\left\{ -\delta\sum_{k=2}^n \left( 1 + \frac{2d-1}{k} \right) \right\}\\
	&\sim \frac{(2d)!}{(n+2)(n+3)...(n+2d)} \exp\left\{ -\delta(n-1+ (2d-1)\log n - (2d-1)) \right\}\\
	&=\frac{(2d)!}{(n+2)(n+3)...(n+2d)} \frac{e^{-\delta n}e^{2\delta d}}{n^{\delta(2d-1)}}
	\end{align*}
	which has the same form as in the one-dimensional case. For $d \geq 2$, the unperturbed walk is diffusive, 
as $\sum_{n=1}^\infty \mP(\tau \geq n) < \infty$. Let $\tau_\delta$ denote the model perturbed 
by $\delta$. By Dominated Convergence Theorem 
\[\lim_{\delta \to 0^+}\mE[\tau_\delta] = \lim_{\delta \to 0^+} 
\sum_{n=1}^\infty \mP(\tau \geq n) = \sum_{n=1}^\infty \lim_{\delta \to 0^+} \mP(\tau \geq n) = \mE[\tau_0]. \]
	Thus $\nu_\delta \to \nu$ as $\delta \to 0^+$.
%, and the model does not produce any residual diffusion. 
For the enhanced diffusivity, by the integral test, it suffices to consider the integral
	\[ \int_1^\infty \frac{e^{-kx}}{{x^{(2d-1)(1+k)}}} dx  =: f(k)\]
	we have
	\[ \frac{\partial}{\partial k}f(k) = \int_1^\infty \frac{e^{-kx}}{x^{(2d-1)k+ 2d-2}}(x + (2d-1)\log x)dx. \]
	Since $d\geq 2$, the integral converges for any non-negative value of $k$. By the Dominated Convergence Theorem,
	\[ \lim_{k \to 0^+} \frac{\partial}{\partial k}f(k)  = 
	\int_1^\infty \lim_{k \to 0^+} \frac{e^{-kx}}{x^{(2d-1)k+ 2d-2}}(x + (2d-1)\log x)dx < \infty, \]
	which implies that $\mE[\tau_\delta]$ grows at a linear rate near $\delta=0$, and therefore
	\[ \nu_\delta = \nu_0 + O(\delta). \]
	%where $M := \mE[\tau_0]$.
	
	\subsubsection{Perturbation under model II and III:}
	Consider the model
	\[
	\mP(S_{n+1}=S_n+x \vert \mF_n) =
	\begin{cases}
	\max\{\frac{1+n}{2d + n}-\delta\xi_n,0\}, \quad \text{if }\{S_n,S_n+x\} = e_n,\\
	\min\{\frac{1}{2d + n}+\delta\xi_n,1\}, \quad \text{if }\{S_n,S_n+x\} \ne e_n.
	\end{cases}
	\]
	If $(\xi_n)_{n \in \mN}$ is a sequence of random variables such that $\mE[\xi_n] < \infty$, for all $n$, 
one can use an analogous argument to that of section \ref{StochasticProof} to show $\nu_\delta \to \nu_0$ at 
the same rate as model I.
% In other words,
%	\[ \nu_\delta = O\left(\frac{1}{M-\delta}\right) \]
	When $\mE[\xi_n] = \infty$, let $f = f_{\xi_n}$. The proof of all three cases are identical. We present the proof of the second case below:	
	\paragraph{Case 2:}$\bm{f(x) = O(x^{-(1+j)})}$, \textbf{for} $\bm{0 < j < 1:}$\\
	Using a similar computation to section \ref{StochasticRate}, we have
	\begin{align*}
	\mP(\tau\geq n) &= \prod_{k=2}^{n}\left( \int_0^{\frac{1+k}{\delta(2d+k)}} \left(\frac{1+k}{2d+k}-\delta x\right)f(x)dx \right)\\
	&\to \prod_{k=2}^n \left(\frac{1+k}{2d+k}- C_1\delta^j \right)\\
	&\to \frac{C_2}{(n+2)(n+3)...(n+2d)} \frac{e^{-\delta^j n}e^{2\delta^j d}}{n^{\delta^j(2d-1)}}
	\end{align*}
	and the Dominated Convergence Theorem guarantees convergence of $\nu_\delta$. For the enhanced diffusivity, it suffices to consider the integral 
	\[ \int_1^\infty \frac{e^{-k^jx}}{{x^{(2d-1)(1+k^j)}}} dx  =: f(k^j)\]
	and 
	\[ \frac{\partial}{\partial k^j}f(k^j) = \int_1^\infty \frac{e^{-k^jx}}{{x^{(2d-1)(1+k^j)}}}(x+(2d-1)\log x) dx < \infty \]
	the integral converges for any non-negative value of $k$. This implies $\mE[\tau_\delta]$ grows at a rate of $\delta^j$ near $\delta=0$. Therefore
	\[ \nu_\delta = \nu_0 + O(\delta^j). \]
	Using an analogous argument, one gets the result for Cases 1 and 3, 
where the construction for Case 3 is the same as in section \ref{StochasticRate}.

	\section{Conclusions}
	The SeRW model in one dimension with identity reinforcement function was found to be diffusive 
when perturbed with a small probability $\delta$ of breaking out of the last traversed edge, 
no matter how small $\delta$ is. The enhanced diffusivity is logarithmically close to  
residual diffusivity as $\delta$ tends to zero.
We studied a few variations of the perturbed models, where the perturbation $\delta\; \xi_n$ is stochastic, and 
    the distribution of $\xi_n$ may or may not depend on $n$. These models intend to create a "fat tail" 
   as $n$ increases so it is more likely for the walk to break out of the last traversed edge. 
  % However none of these models were diffusive enough to sustain residual diffusivity. 
  % As $\delta$ tends to zero, $\nu_\delta$ also converges to zero. 
For most cases, 
the enhanced diffusivity is $\nu_\delta = O\left(\frac{1}{|\log\delta |}\right)$. 
The highest enhanced diffusivity is $\nu_\delta = O\left(\frac{1}{\log|\log\delta|}\right)$.
This was achieved when $\xi_n$ has a very fat tail, $f_{\xi_n}(x) = O\left(\frac{1}{x(\log x)^2}\right)$, which is much fatter than that of the Cauchy distribution.
	In higher dimensions, 
    the baseline SeRW with identity reinforcement function is already diffusive 
    and the enhanced diffusivity reaches a rate as high as  $O(\log^{-2}\delta)$.
	%It is worth noting that in all cases, the perturbed SeRW model achieves a stronger enhanced diffusivity rate than the perturbed simpled random walk model, which has a rate of $O(\delta^2)$.
%\medskip
	
	In future work, we plan to explore dissimilar random walk models with memory mechanism 
    and study enhanced diffusivities.

\section{Acknowledgement}
The authors would like to thank Prof. P. Diaconis for a helpful conversation on reinforced random walk and 
his interest in \cite{LXY_rd}.

\end{document}